\documentclass{amsart}
\usepackage[a4paper,hmargin=30mm,vmargin=30mm]{geometry}
\usepackage{amssymb,mathtools}
\usepackage{enumitem}
\usepackage{bussproofs}
\usepackage[colorlinks=true,linkcolor=blue,citecolor=blue,urlcolor=blue]{hyperref}

\SetEnumitemKey{conds}{leftmargin=3.6em,labelwidth=3.1em,labelsep=0.5em,align=left}

\newtheorem{theorem}{Theorem}[section]
\newtheorem{proposition}[theorem]{Proposition}
\newtheorem{lemma}[theorem]{Lemma}
\newtheorem{corollary}[theorem]{Corollary}
\theoremstyle{definition}
\newtheorem{definition}[theorem]{Definition}
\newtheorem{question}[theorem]{Question}

\newcommand{\RCA}{\mathsf{RCA}_0}
\newcommand{\ACA}{\mathsf{ACA}_0}
\newcommand{\ACAp}{\mathsf{ACA}_0^+}
\newcommand{\ATR}{\mathsf{ATR}_0}
\newcommand{\SAC}{\Sigma^1_1\text{-}\mathsf{AC}_0}
\newcommand{\SDC}{\Sigma^1_1\text{-}\mathsf{DC}_0}
\newcommand{\wRFN}{\omega\mathsf{RFN}}
\newcommand{\WO}{\mathsf{WO}}
\newcommand{\WOP}{\mathsf{WOP}}
\newcommand{\M}{\mathcal{M}}
\newcommand{\N}{\mathfrak{N}}
\newcommand{\U}{\mathfrak{U}}
\newcommand{\Prv}{\omega\mathsf{Prv}}
\newcommand{\Sat}{\mathsf{Sat}}
\newcommand{\pr}{\Box}
\newcommand{\gn}[1]{\ulcorner #1\urcorner}

\begin{document}

\title{Inverting $\omega$-model reflection}
\author{Koshiro Ichikawa}
\address{Department of Mathematical Informatics, Graduate School of Informatics, Nagoya University, Furo-cho, Chikusa-ku, Nagoya 464-8601, Japan}
\email{ichikawa.koshiro.k2@s.mail.nagoya-u.ac.jp}

\begin{abstract}
For a theory $T$, the $\omega$-model reflection principle $\wRFN(T)$ says that every set is contained in a countable coded $\omega$-model of $T$.
We ask which $\Pi^1_2$-sentences are provably equivalent over $\ACA$ to $\wRFN(\ACA+\psi)$ for some sentence $\psi$.
Every such sentence implies $\ACAp$.
Using an analogue of the Friedman--Goldfarb--Harrington theorem for $\omega$-provability, in which witnesses are compared via well-orderings, we show that every $\Pi^1_2$-sentence implying $\ATR$ is of this form, uniformly in the sentence.
On the other hand, theories of the form $\ACA+\wRFN(T)$ are closed under a rule form of $\Pi^1_3$ $\omega$-model reflection for sentences $\forall n\exists Y\varphi$, and hence under a rule of $\Sigma^1_1$-choice.
As a consequence, for every $\Pi^1_2$-sentence $\chi$ that holds in the $\omega$-model of the hyperarithmetical sets, some $\Pi^1_2$-sentence implying $\chi$ is not of this form; this covers $\ACAp$, $\Pi^0_\beta\text{-}\mathsf{CA}_0$ for computable well-orderings $\beta\ge1$, the iterations of $\omega$-model reflection over $\ACA$ along computable well-orderings, and the sentences $\forall X\exists Y\theta(X,Y)$ with $\theta$ arithmetical such that every $X$ has a solution hyperarithmetic in $X$, including the true sentences $\forall X\exists!Y\theta(X,Y)$.
Finally, we show that $\ATR$ is equivalent to $\wRFN(\ACA+\mathsf I)$ for a natural $\Pi^1_2$-sentence $\mathsf I$ built from the Veblen functions.
\end{abstract}

\maketitle

\section{Introduction}

Let $T$ be a consistent computably enumerable extension of $\mathsf{PA}$.
The Friedman--Goldfarb--Harrington (FGH) theorem says that for every $\Sigma_1$-sentence $\sigma$ there is a sentence $\rho$ such that $\mathsf{PA}\vdash\mathrm{Pr}_T(\gn{\rho})\leftrightarrow\sigma\lor\mathrm{Pr}_T(\gn{\bot})$; see Smory\'nski~\cite[p.~366]{Smo81}, Visser~\cite[Theorem~4.1]{Vis84} and Lindstr\"om~\cite[Exercise~2.26(a)]{Lin03}.
Visser~\cite[Section~3]{Vis05} discusses its history, and Kurahashi~\cite{Kur23} and Joosten~\cite{Joo15} study refinements and generalizations.
Taking negations, every $\Pi_1$-sentence $\pi$ with $\mathsf{PA}\vdash\pi\to\mathrm{Con}_T$ is $\mathsf{PA}$-provably equivalent to $\mathrm{Con}(T+\neg\rho)$ for a suitable $\rho$:
the $\Pi_1$-sentences that are consistency statements over $T$ are exactly those implying $\mathrm{Con}_T$.

We study the analogous question for $\omega$-model reflection in second-order arithmetic.
For a finitely axiomatized extension $T$ of $\ACA$, the $\omega$-model reflection principle $\wRFN(T)$ says that every set is contained in a countable coded $\omega$-model of $T$; it is equivalent over $\ACA$ to the $\Pi^1_2$ and the $\Pi^1_1$ $\omega$-model reflection principles for $T$ of Pakhomov and Walsh~\cite[Proposition~3.9 and Theorem~3.10]{PW23}.
Many familiar principles have this form.
For example, $\ACAp$ is equivalent to $\wRFN(\ACA)$~\cite[Proposition~3.9]{PW23} and to the well-ordering principle $\forall X(\WO(X)\to\WO(\varepsilon_X))$~\cite{AR09,MM11}, and Rathjen~\cite[Theorem~1.4]{Rat14} shows that the well-ordering principle $\forall X(\WO(X)\to\WO(\Gamma_X))$ is equivalent to the statement that every set is contained in a countable coded $\omega$-model of $\mathsf{ATR}$.
Similarly, for a theory that is not finitely axiomatized, Suzuki~\cite{Suz25} shows that, for $n>1$, the $\Sigma^0_n$ relative leftmost path principle of Towsner is equivalent over $\RCA$ to the $\omega$-model reflection of $\Pi^1_{n+1}$ transfinite induction.

We ask which $\Pi^1_2$-sentences $\varphi$ satisfy $\ACA\vdash\varphi\leftrightarrow\wRFN(\ACA+\psi)$ for some sentence $\psi$.
Such a $\varphi$ implies $\wRFN(\ACA)$, that is, $\ACAp$.
In the arithmetical setting the corresponding necessary condition is also sufficient.
Our first result is a partial analogue.

\begin{theorem}\label{thm:A}
There is a primitive recursive function $\varphi\mapsto\eta_\varphi$ from $\Pi^1_2$-sentences to $\Pi^1_2$-sentences such that $\ACA\vdash\wRFN(\ACA+\eta_\varphi)\to\varphi$ and $\ATR\vdash\varphi\to\wRFN(\ACA+\eta_\varphi)$ for every $\Pi^1_2$-sentence $\varphi$.
In particular, every $\Pi^1_2$-sentence $\varphi$ with $\ACA\vdash \varphi\to\ATR$ is $\ACA$-provably equivalent to $\wRFN(\ACA+\eta_\varphi)$.
\end{theorem}

For a $\Pi^1_2$-sentence $\chi$, we call the class of $\Pi^1_2$-sentences $\varphi$ with $\ACA\vdash\varphi\to\chi$ the \emph{cone above $\chi$}, and we say that $\omega$-model reflection is \emph{invertible on the cone above $\chi$} if every sentence in this cone is $\ACA$-provably equivalent to $\wRFN(\ACA+\psi)$ for some sentence $\psi$ (Definition~\ref{def:main}).
The cone above $\chi$ is the set of $\Pi^1_2$-sentences in the cone generated by $\chi$ in the sense of Walsh~\cite[Definition~2.3]{Wal19}, and its elements need not be true.
Thus Theorem~\ref{thm:A} says that $\omega$-model reflection is invertible on the cone above $\ATR$, which is a true cone in Walsh's sense.

The proof (Section~\ref{sec:fgh}) follows the FGH theorem with $\mathrm{Pr}_T$ replaced by $\omega$-provability, the $\Sigma^1_2$ dual of $\omega$-model reflection.
A $\Sigma^1_2$-sentence can be written as $\exists X\WO(L_X)$, and we compare witnesses by comparing the well-orderings $L_X$.
In the FGH theorem the comparison is shown to be exclusive inside the theory, which here would require comparability of well-orderings in the $\omega$-models of the reflected theory.
Instead, we use $\ATR$ only outside the $\omega$-models.
When reflection fails at a set $Y$, a genuine well-ordering $K_Y$ is determined, and by a result of Avigad and Sommer~\cite{AS99}, $\ATR$ provides an $\omega$-model of $\ACA$ in which every well-ordering is comparable with $K_Y$.

Our second result, Theorem~\ref{thm:B} below, is a closure property of the theories $\ACA+\wRFN(T)$, from which we derive that $\omega$-model reflection is not invertible on the cone above any $\Pi^1_2$-sentence that holds in the $\omega$-model $\mathrm{HYP}$ of the hyperarithmetical sets.
These cones include those above $\ACAp$, above the systems $\Pi^0_\beta\text{-}\mathsf{CA}_0$ of Marcone and Montalb\'an~\cite{MM11} for computable well-orderings $\beta\ge1$, which include the iterations of $\omega$-model reflection over $\ACA$ along computable well-orderings, and above any sentence $\forall X\exists Y\theta(X,Y)$ with $\theta$ arithmetical such that every $X$ has a solution hyperarithmetic in $X$, which includes every true sentence $\forall X\exists!Y\theta(X,Y)$.

\begin{theorem}\label{thm:B}
Let $T$ be a finitely axiomatized extension of $\ACA$, put $S:=\ACA+\wRFN(T)$, and let $\varphi(n,Y,Z)$ be a $\Sigma^1_1$-formula.
Then $S$ is closed under the rule
\begin{prooftree}
\AxiomC{$\forall Z\forall n\exists Y\varphi(n,Y,Z)$}
\UnaryInfC{$\forall Z\exists\M\bigl(Z\in\M\land\M\models T+\forall n\exists Y\varphi(n,Y,Z)\bigr)$}
\end{prooftree}
and hence under the rule with the same premise and the conclusion $\forall Z\exists W\forall n\varphi(n,(W)_n,Z)$.
\end{theorem}

Theorem~\ref{thm:B}, proved in Section~\ref{sec:choice}, is a rule form of a reflection principle of Pakhomov and Walsh.
Their principle $\Pi^1_3\text{-}\omega\mathsf{RFN}(T)$ says that every $\Pi^1_3$-sentence with set parameters that holds in all $\omega$-models of $T$ is true; equivalently, every true $\Sigma^1_3$-sentence with set parameters holds in some $\omega$-model of $T$~\cite[Section~2.2]{PW23}.
Applied to the true $\Pi^1_2$-sentence $\forall n\exists Y\varphi(n,Y,Z)$, it gives an $\omega$-model as in the conclusion of Theorem~\ref{thm:B}, and hence the corresponding instance of $\SAC$.
The theory $S$ need not prove these instances (for $T=\ACA$, see Proposition~\ref{prop:IF} and the proof of Corollary~\ref{cor:examples} with $\beta=\omega$), but Theorem~\ref{thm:B} shows that it is closed under the corresponding rule.
The restriction to sentences of the form $\forall n\exists Y\varphi$ cannot be dropped; see the remark after the proof of Theorem~\ref{thm:B}.
Rules of this kind appear in the variant~\cite[Lemma~2.11]{PW23} of Beklemishev's reduction property~\cite[Lemma~4.2]{Bek03}: Pakhomov and Walsh show that the closure of a $\Pi^1_2$-axiomatized theory under the rule
\begin{prooftree}
\AxiomC{$\psi$}
\UnaryInfC{$\wRFN(T+\psi)$}
\end{prooftree}
for $\Pi^1_2$-sentences $\psi$ (possibly with set constants) has the same $\Pi^1_2$-consequences as its extension by $\Pi^1_3\text{-}\omega\mathsf{RFN}(T)$.
Theorem~\ref{thm:B} says that $S$ itself is closed under this rule for $\psi$ of the form $\forall n\exists Y\varphi$.
Together with~\cite[Lemma~2.10]{PW23}, its proof also shows that $S+\SAC$ is $\Pi^1_2$-conservative over $S$ (Corollary~\ref{cor:cons}), an analogue of the theorem of Barwise and Schlipf~\cite{BS75} that $\SAC$ is $\Pi^1_2$-conservative over $\ACA$ (see also~\cite[Theorem~IX.4.4]{Sim09}).
If $\omega$-model reflection is invertible on the cone above $\chi$, then $\ACA+\chi$ proves the instances of $\SAC$ without set parameters (Corollary~\ref{cor:necessary}).
Every $\omega$-model of these instances contains all hyperarithmetical sets, whereas if $\chi$ holds in $\mathrm{HYP}$, then $\ACA+\chi$ has a hyperarithmetical $\omega$-model, which does not contain all hyperarithmetical sets; so $\chi$ does not hold in $\mathrm{HYP}$ (Proposition~\ref{prop:HYP}).
Since $\ATR$ proves $\SAC$, these necessary conditions are consistent with Theorem~\ref{thm:A}.
Applying Theorem~\ref{thm:B} to finite iterations of a suitable arithmetically defined operation, we also show that natural sentences in some of these cones, such as those above $\Pi^0_\beta\text{-}\mathsf{CA}_0$, are not of the form $\wRFN(\ACA+\psi)$ (Proposition~\ref{prop:IF} and Corollary~\ref{cor:examples}).

Hashimoto~\cite{Has24} studied the case $\varphi=\ATR$, regarded as a $\Pi^1_2$-sentence (see Section~\ref{sec:prelim}).
He writes $\mathsf{RFN}^{-1}(\ATR)$ for the class of theories $T$ such that the statement that every set is contained in a countable coded $\omega$-model of $T+\ACA$ is equivalent to $\ATR$ over $\ACA$~\cite[Definition~5.6]{Has24}.
He shows that this class contains every theory between $\mathsf{JI}_0$ and $\SDC$~\cite[Proposition~5.8]{Has24}, and that it contains no sentence of the form $\forall X\exists!Y\theta(X,Y)$ with $\theta$ arithmetical~\cite[Corollary~5.13]{Has24}.
Since no theory of hyperarithmetic analysis is finitely axiomatized by $\Sigma^1_3$-sentences~\cite[Theorem~5.5]{Has24}, he asks whether $\mathsf{RFN}^{-1}(\ATR)$ contains a $\Sigma^1_3$-sentence~\cite[Question~6.1]{Has24}.
Theorem~\ref{thm:A} with $\varphi=\ATR$ gives a $\Pi^1_2$-sentence $\eta_{\ATR}$, in particular a $\Sigma^1_3$-sentence, in $\mathsf{RFN}^{-1}(\ATR)$, and this answers Hashimoto's question positively.
The sentence $\eta_{\ATR}$ is obtained by diagonalization.
In Section~\ref{sec:atr} we give a natural one: $\ATR$ is equivalent to $\wRFN(\ACA+\mathsf I)$, where $\mathsf I$ says that $\varphi_X$ preserves well-orderings whenever $\varphi_X(0)$ is well-ordered (Theorem~\ref{thm:atr}).
The proof reduces $\wRFN(\ACA+\mathsf I)$ to iterated $\Pi^1_1$-reflection by a theorem of Pakhomov and Walsh~\cite{PW23}, and bounds the $\Pi^1_1$-strength of the iterates using their computation of the proof-theoretic ordinals of iterated $\omega$-model reflection together with well-ordering assumptions.

\subsection*{Historical remark}
The constructions of Theorem~\ref{thm:A} and Section~\ref{sec:atr} were found at the Oberwolfach workshop mentioned in the acknowledgements.
On 15 November 2026, Yudai Suzuki told the author that~\cite[Question~6.1]{Has24} had been answered positively; the answer had been obtained when he asked ChatGPT~6 Astra (OpenAI) about it.
On 16 November, the author found the sentence $\mathsf I$ of Section~\ref{sec:atr}.
On 17 November, the diagonal construction of Theorem~\ref{thm:A} was found independently by the author and by Keita Yokoyama, who had heard of the sentence $\mathsf I$.

\section{Preliminaries}\label{sec:prelim}

We work in the language $L_2$ of second-order arithmetic and follow Simpson~\cite{Sim09} for subsystems and their basic properties.
A \emph{(countable coded) $\omega$-model} is a set $\M$, viewed as the sequence of sets $(\M)_j$; we write $X\in\M$ if $X=(\M)_j$ for some $j$.
For a formula $\varphi$ with set parameters in $\M$, $\M\models\varphi$ denotes the relativization of $\varphi$ to $\M$; number quantifiers are not relativized.
This is a separate formula for each $\varphi$: $\ACA$ does not prove that every $\omega$-model has a full satisfaction class, which is available in $\ACAp$~\cite[Section~2.2]{PW23}.
Thus our $\omega$-models are the weak $\omega$-models of~\cite[Section~3.3]{PW23}, whereas the $\omega$-models of~\cite{PW23} carry full satisfaction classes.
Since $\ACA$ is finitely axiomatizable (see e.g.\ \cite[proof of Lemma~4.1]{Rat14}), for a finitely axiomatized extension $T=\ACA+\psi$ the statement $\M\models T$ is arithmetical in $\M$.
All theories in this paper are finitely axiomatized extensions of $\ACA$, except the iterated reflection theories $\Pi^1_2\text{-}\omega\mathbf R^\alpha(\ACA)$ of~\cite{PW23} and theories in results quoted from the literature; for these, $\M\models T$ and $\wRFN(T)$ are understood as in~\cite[Section~2.2]{PW23}.

We write $\gn{\varphi}$ for the code of a formula $\varphi$.
In $\ACA$ there are partial truth predicates for $\Pi^1_n$-formulas with set parameters~\cite[Section~2.2]{PW23}; let $\mathsf{Tr}_{\Pi^1_3}(x)$ be the one for $\Pi^1_3$-sentences, that is, with no parameters.
It is a $\Pi^1_3$-formula, and $\ACA\vdash\varphi\leftrightarrow\mathsf{Tr}_{\Pi^1_3}(\gn{\varphi})$ for every $\Pi^1_3$-sentence $\varphi$.
Let $\Sat(\M,x)$ be the relativization of $\mathsf{Tr}_{\Pi^1_3}(x)$ to $\M$.
It is arithmetical in $\M$ and $x$, and relativizing the proof of $\varphi\leftrightarrow\mathsf{Tr}_{\Pi^1_3}(\gn{\varphi})$ to $\M$ shows that, for every $\Pi^1_3$-sentence $\varphi$, $\ACA$ proves $\Sat(\M,\gn{\varphi})\leftrightarrow\M\models\varphi$ for every $\omega$-model $\M$ of $\ACA$; no full satisfaction class is needed.
We apply $\Sat$ only to codes of sentences, namely of Boolean combinations of $\Sigma^1_2$-sentences, each of which we identify with a logically equivalent $\Pi^1_3$-sentence obtained primitive recursively from it.

\begin{definition}\label{def:main}
Let $T$ be a finitely axiomatized extension of $\ACA$.
\begin{enumerate}[label=(\roman*)]
\item $\wRFN(T)$ is the sentence $\forall X\exists\M(X\in\M\land\M\models T)$.
\item $\Prv(T,x)$ is the formula $\exists X\forall\M(X\in\M\land\M\models T\to\Sat(\M,x))$.
\item For a $\Pi^1_2$-sentence $\chi$, the \emph{cone above $\chi$} is the class of $\Pi^1_2$-sentences $\varphi$ with $\ACA\vdash\varphi\to\chi$, and \emph{$\omega$-model reflection is invertible on the cone above $\chi$} if every sentence in this cone is $\ACA$-provably equivalent to $\wRFN(\ACA+\psi)$ for some sentence $\psi$.
\end{enumerate}
\end{definition}

Both $\wRFN(T)$ and $\Prv(T,x)$ are formulated in $\ACA$; the former is $\Pi^1_2$ and the latter $\Sigma^1_2$.
By the property of $\Sat$, $\ACA\vdash\Prv(T,\gn{\varphi})\leftrightarrow\exists X\forall\M(X\in\M\land\M\models T\to\M\models\varphi)$ for every $\Pi^1_3$-sentence $\varphi$.
For a sentence $\varphi$ we write $\pr_T\varphi$ for $\Prv(T,\gn{\varphi})$.
Over $\ACA$, $\wRFN(T)$ implies $\ACAp$ by~\cite[Proposition~3.9]{PW23}, and in $\ACAp$ every $\omega$-model has a full satisfaction class, which agrees with relativization on each fixed formula.
Hence $\wRFN(T)$ is $\ACA$-provably equivalent to the statement that every set is contained in an $\omega$-model of $T$ in the sense of~\cite{PW23}, and so, by~\cite[Theorem~3.10]{PW23}, to the $\Pi^1_2$ $\omega$-model reflection principle $\Pi^1_2\text{-}\omega\mathsf{RFN}(T)$ of~\cite[Section~2.2]{PW23}.
For every Boolean combination $\psi$ of $\Sigma^1_2$-sentences, by the property of $\Sat$,
\begin{equation}\label{eq:dual}
\ACA\vdash\neg\pr_T\neg\psi\leftrightarrow\wRFN(T+\psi),\qquad
\ACA\vdash\neg\pr_T\bot\leftrightarrow\wRFN(T).
\end{equation}
Dually, for every Boolean combination $\varphi$ of $\Sigma^1_2$-sentences, $\pr_T\varphi$ is $\ACA$-provably equivalent to the $\Sigma^1_2$ $\omega$-provability predicate $\Sigma^1_2\text{-}\omega\mathsf{Prv}(T,\varphi)$ of~\cite[Section~2.2]{PW23}, which is defined as the negation of $\Sigma^1_2\text{-}\omega\mathsf{Con}(T+\neg\varphi)$ and says that some true $\Sigma^1_2$-sentence with set parameters holds in no $\omega$-model of $T+\neg\varphi$.
Indeed, $\Sigma^1_2\text{-}\omega\mathsf{Con}(T+\neg\varphi)$ is $\ACA$-provably equivalent to the $\Pi^1_2$ $\omega$-model reflection principle for $T+\neg\varphi$ in the sense of~\cite[Section~2.2]{PW23}, hence to $\wRFN(T+\neg\varphi)$, whose negation is equivalent to $\pr_T\varphi$ by~\eqref{eq:dual}.
Since $\Prv(T,x)$ is a formula of the code $x$, the diagonal lemma applies to it in Section~\ref{sec:fgh}.

\begin{proposition}\label{prop:basic}
Let $T$ be a finitely axiomatized extension of $\ACA$, and let $\varphi,\psi$ range over Boolean combinations of $\Sigma^1_2$-sentences.
\begin{enumerate}[label=(\roman*)]
\item If $T\vdash\varphi$, then $\ACA\vdash\pr_T\varphi$.
\item $\ACA\vdash\pr_T(\varphi\to\psi)\to(\pr_T\varphi\to\pr_T\psi)$.
\item $\ACA\vdash\sigma\to\pr_T\sigma$ for every $\Sigma^1_2$-sentence $\sigma$.
\end{enumerate}
\end{proposition}

The proofs are straightforward.
Applied to the $\Sigma^1_2$-sentence $\pr_T\varphi$, item~(iii) gives $\ACA\vdash\pr_T\varphi\to\pr_T\pr_T\varphi$, so (i)--(iii) are the Hilbert--Bernays--L\"ob derivability conditions for $\pr_T$ over $\ACA$.
Hence, by the usual argument (see~\cite{Boo93}), L\"ob's theorem holds for $\pr_T$, and the provability logic $\mathsf{GL}$ is sound for it: if $\mathsf{GL}\vdash A$, then $\ACA\vdash A^*$ for every interpretation $^*$ that maps propositional variables to Boolean combinations of $\Sigma^1_2$-sentences, commutes with the Boolean connectives and maps $\Box B$ to $\pr_T B^*$.

For linear orders $L$ and $K$, we write $K\lhd L$ if $K$ is isomorphic to a proper initial segment of $L$, and $K\unlhd L$ if $K$ is isomorphic to an initial segment of $L$.
Both relations are $\Sigma^1_1$.
If $\WO(L)$ and $K\lhd L$ or $K\unlhd L$, then $\WO(K)$, provably in $\RCA$, since an isomorphism maps descending sequences to descending sequences.
By the Kleene normal form theorem and the Kleene--Brouwer ordering (see \cite[Section~V.1]{Sim09}), there is a primitive recursive function $e$ such that every $\Sigma^1_2$-sentence $\sigma$ is $\ACA$-provably equivalent to $\exists X\WO(L_{e(\gn{\sigma}),X})$, where $L_{e,X}$ is the Kleene--Brouwer ordering of a tree computable in $X$ uniformly in $e$ and $X$.
We call a sentence $\exists X\WO(L_{e,X})$ that is $\ACA$-provably equivalent to $\sigma$ a \emph{Kleene normal form} of $\sigma$.

For a linear order $\beta$ with a first element and sets $X,Y$, we write $X=Y^{(\beta)}$ if $X$ is the jump hierarchy along $\beta$ starting with $Y$, as defined in~\cite[Section~2]{MM11}; this relation is arithmetical in $X$, $Y$ and $\beta$.
For a computable well-ordering $\beta\ge1$, the system $\Pi^0_\beta\text{-}\mathsf{CA}_0$ is $\RCA$ together with the statement that $\beta$ is well-ordered and $\forall Y\exists X(X=Y^{(\beta)})$~\cite[Section~2]{MM11}.
Thus $\Pi^0_1\text{-}\mathsf{CA}_0$ is $\ACA$ and $\Pi^0_\omega\text{-}\mathsf{CA}_0$ is $\ACAp$.
The theory $\ATR$ is axiomatized over $\RCA$ by the single $\Pi^1_2$-sentence $\forall\alpha(\WO(\alpha)\to\forall Y\exists X(X=Y^{(\alpha)}))$~\cite[Section~2]{MM11}.
As in~\cite[Section~2]{PW23}, we identify these finitely axiomatized theories with the sentences axiomatizing them.
By a theorem of Friedman (see~\cite[Theorem~V.6.8]{Sim09} and~\cite{FH90}), $\ATR$ proves that any two well-orderings $L,K$ are comparable; in particular $K\lhd L$ or $L\unlhd K$.

We use $\ATR$ in Section~\ref{sec:fgh} only through the following proposition.
Here $\M\models L\unlhd K$ says that $\M$ contains an isomorphism from $L$ onto an initial segment of $K$, and similarly for $\M\models K\lhd L$.

\begin{proposition}\label{lem:cmp}
$\ATR$ proves
\[
\forall K\forall Z\Bigl[\WO(K)\to\exists\M\bigl(K,Z\in\M\land\M\models\ACA\land\forall L\in\M\bigl(\M\models\WO(L)\to\M\models L\unlhd K\lor K\lhd L\bigr)\bigr)\Bigr].
\]
\end{proposition}

\begin{proof}
By~\cite[Lemma~7.5]{AS99} (see also~\cite[Lemma~VIII.4.19]{Sim09}), $\ATR$ proves that every set is contained in a countable coded $\omega$-model of $\SAC$; let $\M$ be one containing $K\oplus Z$.
Then $\M$ contains the solution of every arithmetical transfinite recursion along $K$ with parameters in $\M$.
Indeed, if $a$ were the least element of $K$ such that the recursion along $K{\restriction}a$ has no solution in $\M$, then the solutions along $K{\restriction}b$ for $b<_K a$ would lie in $\M$, $\SAC$ in $\M$ would collect them, and arithmetical comprehension in $\M$ would give a solution along $K{\restriction}a$.
In the same way, $\SAC$ in $\M$ collects the solutions along all $K{\restriction}a$, and arithmetical comprehension in $\M$ gives a solution along $K$.

Now let $L\in\M$ with $\M\models\WO(L)$, and define $f$ by arithmetical transfinite recursion along $K$: $f(a)$ is the $L$-least element of $L\setminus f[K{\restriction}a]$, and $f(a)$ is undefined if this set is empty.
Then $f\in\M$.
Each set $L\setminus f[K{\restriction}a]$ belongs to $\M$, so it has a least element whenever it is nonempty, since $\M\models\WO(L)$.
If $L\setminus f[K{\restriction}a]$ is empty for some $a$, then for the least such $a$ the map $f$ is an isomorphism from $K{\restriction}a$ onto $L$, so its inverse witnesses $\M\models L\unlhd K$.
Otherwise $f$ is an isomorphism from $K$ onto an initial segment of $L$; if this segment is $L$, then the inverse of $f$ witnesses $\M\models L\unlhd K$, and otherwise $f$ witnesses $\M\models K\lhd L$.
\end{proof}

\section{\texorpdfstring{An FGH theorem for $\omega$-provability}{An FGH theorem for omega-provability}}\label{sec:fgh}

\begin{definition}
For $\Sigma^1_2$-sentences in Kleene normal form $\sigma\equiv\exists X\WO(L_X)$ and $\tau\equiv\exists Y\WO(K_Y)$, let
\begin{align*}
\sigma\preccurlyeq\tau&:\equiv\exists X\bigl[\WO(L_X)\land\neg\exists Y(K_Y\lhd L_X)\bigr],\\
\sigma\prec\tau&:\equiv\exists X\bigl[\WO(L_X)\land\neg\exists Y(K_Y\unlhd L_X)\bigr].
\end{align*}
\end{definition}

This mirrors the witness comparison of $\Sigma_1$-sentences (see~\cite[Lemma~1.3]{Lin03} and~\cite[Definition~1.3]{Kur23}); here witnesses are measured by the order types of $L_X$ and $K_Y$: $\sigma\preccurlyeq\tau$ says that some witness $X$ of $\sigma$ has no witness $Y$ of $\tau$ with $K_Y$ shorter than $L_X$, and $\sigma\prec\tau$ says that it has none with $K_Y$ at most as long as $L_X$.
The inner subformulas do not mention $\WO(K_Y)$, which follows from $\WO(L_X)$; this ensures that $\sigma\preccurlyeq\tau$ and $\sigma\prec\tau$ are $\Sigma^1_2$-sentences.

\begin{proposition}\label{prop:wc}
Let $\sigma\equiv\exists X\WO(L_X)$ and $\tau\equiv\exists Y\WO(K_Y)$ be $\Sigma^1_2$-sentences in Kleene normal form.
\begin{enumerate}[label=(\roman*)]
\item $\ACA\vdash\sigma\preccurlyeq\tau\to\sigma$.
\item $\ACA\vdash\sigma\prec\tau\to\sigma\preccurlyeq\tau$.
\item $\ACA\vdash\sigma\land\neg\tau\to\sigma\prec\tau$.
\item $\ATR\vdash\neg(\sigma\prec\tau\land\tau\preccurlyeq\sigma)$.
\item $\ATR\vdash\forall Y\bigl[\WO(K_Y)\to\exists\M\bigl(Y\in\M\land\M\models\ACA\land(\M\models\sigma\preccurlyeq\tau\to\sigma)\bigr)\bigr]$.
\end{enumerate}
\end{proposition}

\begin{proof}
(i) is immediate from the definition, and (ii) holds since $K_Y\lhd L_X$ implies $K_Y\unlhd L_X$.
(iii) Let $\WO(L_X)$.
If $K_Y\unlhd L_X$ for some $Y$, then $\WO(K_Y)$ and hence $\tau$.
So under $\neg\tau$, the set $X$ witnesses $\sigma\prec\tau$.

(iv) We argue in $\ATR$.
Suppose that $X$ witnesses $\sigma\prec\tau$ and $Y$ witnesses $\tau\preccurlyeq\sigma$.
Since $L_X$ and $K_Y$ are well-orderings, $K_Y\unlhd L_X$ or $L_X\lhd K_Y$ by the comparability of well-orderings~\cite[Theorem~V.6.8]{Sim09}.
The first contradicts the choice of $X$, and the second that of $Y$.

(v) We argue in $\ATR$.
Let $\WO(K_Y)$, and let $\M$ be as in Proposition~\ref{lem:cmp} for $K=K_Y$ and $Z=Y$.
Suppose $\M\models\sigma\preccurlyeq\tau$.
Then there is $X\in\M$ such that $\M\models\WO(L_X)$ and, since $Y\in\M$, $\M\models\neg(K_Y\lhd L_X)$.
By Proposition~\ref{lem:cmp}, $\M\models L_X\unlhd K_Y$, and hence $L_X\unlhd K_Y$.
Since $K_Y$ is well-ordered, so is $L_X$, and hence $\sigma$.
\end{proof}

Items (i)--(iv) are the analogues of the properties of the witness comparison of $\Sigma_1$-sentences used in the FGH theorem.
There, exclusivity (iv) is used inside the theory, which here would require (iv) in $\omega$-models of $\ACA$, where comparability of well-orderings may fail.
Item~(v) replaces it: once a witness $Y$ of $\tau$ is fixed, $K_Y$ is a genuine well-ordering, and in a suitable $\omega$-model $\sigma\preccurlyeq\tau$ forces $\sigma$.

We now fix a $\Pi^1_2$-sentence $\varphi$, and let $\sigma\equiv\exists X\WO(L_X)$ be a Kleene normal form of $\neg\varphi$.
We write $\pr$ for $\pr_{\ACA}$.
The formula $\forall\M(Y\in\M\land\M\models\ACA\to\Sat(\M,x))$ is $\Pi^1_1$, so by the Kleene normal form theorem there is a linear order $K_{x,Y}$, computable in $Y$ uniformly in $x$ and $Y$, such that $\ACA\vdash\forall x\forall Y\bigl[\WO(K_{x,Y})\leftrightarrow\forall\M(Y\in\M\land\M\models\ACA\to\Sat(\M,x))\bigr]$.
Thus, for every Boolean combination $\rho$ of $\Sigma^1_2$-sentences, $\exists Y\WO(K_{\gn{\rho},Y})$ is a Kleene normal form of $\pr\rho$, and we use it in $\sigma\preccurlyeq\pr\rho$.
By the diagonal lemma applied to the formula $\exists X\bigl[\WO(L_X)\land\neg\exists Y(K_{x,Y}\lhd L_X)\bigr]$, there is a $\Sigma^1_2$-sentence $\rho=\rho_\varphi$, obtained primitive recursively from $\varphi$, such that
\begin{equation}\label{eq:fixed}
\ACA\vdash\rho\leftrightarrow(\sigma\preccurlyeq\pr\rho).
\end{equation}
For this $\rho$ we write $K_Y$ for $K_{\gn{\rho},Y}$.

\begin{theorem}\label{thm:fgh}
Let $\varphi$ and $\rho$ be as above.
\begin{enumerate}[label=(\roman*)]
\item $\ACA\vdash\neg\varphi\to\pr\rho$.
\item $\ATR\vdash\pr\rho\to\neg\varphi$.
\end{enumerate}
\end{theorem}

\begin{proof}
(i) We argue in $\ACA$ and assume $\sigma$.
If $\neg\pr\rho$, then $\sigma\preccurlyeq\pr\rho$ by Proposition~\ref{prop:wc}(iii) and~(ii), so $\rho$ by~\eqref{eq:fixed}, and $\pr\rho$ by Proposition~\ref{prop:basic}(iii).
Hence $\pr\rho$ in either case.

(ii) We argue in $\ATR$ and assume $\pr\rho$.
Let $Y$ witness $\pr\rho$; then $\WO(K_Y)$ by the choice of $K_Y$.
Let $\M$ be as in Proposition~\ref{prop:wc}(v) for $\tau=\pr\rho$.
Since $Y$ witnesses $\pr\rho$ and $\M$ is an $\omega$-model of $\ACA$ containing $Y$, we have $\Sat(\M,\gn{\rho})$ and hence $\M\models\rho$, and since $\ACA$ proves~\eqref{eq:fixed}, $\M\models\sigma\preccurlyeq\pr\rho$.
Hence $\sigma$ by Proposition~\ref{prop:wc}(v), that is, $\neg\varphi$.
\end{proof}

\begin{proof}[Proof of Theorem~\ref{thm:A}]
Put $\eta_\varphi:\equiv\neg\rho_\varphi$, which is $\Pi^1_2$ and obtained primitive recursively from $\varphi$.
By~\eqref{eq:dual} and Proposition~\ref{prop:basic}, $\ACA\vdash\neg\pr\rho_\varphi\leftrightarrow\wRFN(\ACA+\eta_\varphi)$, so Theorem~\ref{thm:fgh}(i) gives the first claim and Theorem~\ref{thm:fgh}(ii) the second.
If moreover $\ACA\vdash\varphi\to\ATR$, then $\ACA+\varphi$ proves $\ATR$, hence $\wRFN(\ACA+\eta_\varphi)$.
Thus $\omega$-model reflection is invertible on the cone above $\ATR$.
\end{proof}

The proof is uniform in parameters, as in the parameterized version of the FGH theorem (see~\cite[Section~1]{Kur23}).
Let $\varphi(v,X)$ be a $\Pi^1_2$-formula with a number variable $v$ and a set variable $X$, and write $\wRFN(\ACA+\eta(\dot v,\dot X))$ for $\forall Y\exists\M\bigl(X,Y\in\M\land\M\models\ACA+\eta(v,X)\bigr)$.
Running the argument above with a partial truth predicate with set parameters in place of $\mathsf{Tr}_{\Pi^1_3}$ and with the diagonal lemma with parameters gives a $\Pi^1_2$-formula $\eta(v,X)$ such that $\ACA\vdash\forall v\forall X\bigl(\wRFN(\ACA+\eta(\dot v,\dot X))\to\varphi(v,X)\bigr)$ and $\ATR\vdash\forall v\forall X\bigl(\varphi(v,X)\to\wRFN(\ACA+\eta(\dot v,\dot X))\bigr)$; we omit the details.
Applied to a partial truth predicate for $\Pi^1_2$-formulas with set parameters, this gives a single formula $\eta$ that works for all $\Pi^1_2$-formulas.

\section{\texorpdfstring{A rule form of $\Pi^1_3$ $\omega$-model reflection}{A rule form of Pi13 omega-model reflection}}\label{sec:choice}

We prove Theorem~\ref{thm:B} and derive from it necessary conditions for invertibility.

\begin{proof}[Proof of Theorem~\ref{thm:B}]
Since $\varphi$ is $\Sigma^1_1$, we have $\varphi(n,Y,Z)\equiv\exists V\varphi'(n,Y,V,Z)$ with $\varphi'$ arithmetical.
By pairing, the two set quantifiers of $\exists Y\exists V\varphi'$ can be merged: $\exists Y\varphi(n,Y,Z)$ is $\ACA$-provably equivalent to $\exists W\varphi'(n,(W)_0,(W)_1,Z)$, and the same holds in every $\omega$-model of $\ACA$, since such models are closed under $\oplus$ and under $W\mapsto(W)_i$.
So both rules for $\varphi$ follow from those for the arithmetical formula $\varphi'(n,(W)_0,(W)_1,Z)$, and we may assume that $\varphi$ is arithmetical.
Suppose that $S\vdash\forall Z\forall n\exists Y\varphi(n,Y,Z)$ and $S\nvdash\forall Z\exists\M(Z\in\M\land\M\models T+\forall n\exists Y\varphi(n,Y,Z))$.
Then there are a model $\U$ of $S$ and a set $Z_0$ of $\U$ such that, in $\U$, no $\omega$-model of $T$ containing $Z_0$ satisfies $\forall n\exists Y\varphi(n,Y,Z_0)$.

Fix a standard number $k$.
Applying $\wRFN(T)$ in $\U$ repeatedly, we obtain, in $\U$, $\omega$-models $\M_0,\dots,\M_k$ of $T$ with $Z_0\in\M_0$ and $\M_i\in\M_{i+1}$ for $i<k$.
Since $\M_{i+1}$ is closed under relative computability, every set of $\M_i$ belongs to $\M_{i+1}$, so $\M_0\subseteq\dots\subseteq\M_k$.
By the choice of $Z_0$, there is a number $b$ in $\U$ such that $\neg\varphi(b,(\M_k)_j,Z_0)$ for all $j$, and hence $\neg\varphi(b,(\M_i)_j,Z_0)$ for all $i\le k$ and all $j$.

We add new set constants $\underline{Z},\underline{\M_0},\underline{\M_1},\dots$ and a new number constant $\underline{b}$.
Let $\Gamma$ consist of $S$ together with the sentences $\underline{Z}\in\underline{\M_0}$, $\underline{\M_i}\models T$, $\underline{\M_i}\in\underline{\M_{i+1}}$ and $\forall j\neg\varphi(\underline{b},(\underline{\M_i})_j,\underline{Z})$ for all $i\in\omega$, where $\underline{\M_i}\models T$ abbreviates the arithmetical formula expressing that $\underline{\M_i}$ is an $\omega$-model of $T$.
By the previous paragraph, with $\underline{Z}$, $\underline{\M_i}$ and $\underline{b}$ interpreted as $Z_0$, $\M_i$ and $b$, every finite subset of $\Gamma$ is satisfiable, so $\Gamma$ has a model $\U^*$ by compactness.
Let $\N$ be the structure with the number part of $\U^*$ whose sets are the sets $(\underline{\M_i})_j$ of $\U^*$, for $i\in\omega$ and $j$ ranging over the numbers of $\U^*$.

We show that $\N$ is a model of $S$.
Arithmetical formulas have the same truth value in $\N$ and in $\U^*$.
Any finitely many sets of $\N$ lie in a single $\underline{\M_i}$.
Since $\underline{\M_i}\models\ACA$ holds in $\U^*$, every set arithmetically definable from them is some $(\underline{\M_i})_j$, so $\N$ satisfies arithmetical comprehension; the induction axiom holds in $\N$ because it holds in $\U^*$.
Hence $\N\models\ACA$.
For $X\in\underline{\M_i}$, the set $\underline{\M_i}$ is itself a set of $\N$ because $\underline{\M_i}\in\underline{\M_{i+1}}$, and the arithmetical statement that $\underline{\M_i}$ is an $\omega$-model of $T$ containing $X$ holds in $\N$.
Hence $\N\models\wRFN(T)$.

Since $\N\models S$, we get $\N\models\exists Y\varphi(\underline{b},Y,\underline{Z})$.
But every set of $\N$ is of the form $(\underline{\M_i})_j$, and $\neg\varphi(\underline{b},(\underline{\M_i})_j,\underline{Z})$ holds in $\U^*$ and thus in $\N$, a contradiction.

For the second rule, let $\M$ be as in the conclusion of the first.
The function $g(n):=\min\{j:\varphi(n,(\M)_j,Z)\}$ is arithmetical in $\M\oplus Z$, so the set $W$ with $(W)_n=(\M)_{g(n)}$ exists by arithmetical comprehension.
\end{proof}

The proof uses that $\forall n\exists Y\varphi(n,Y,Z_0)$ fails in $\M_k$ at a number $b$: since $\M_0\subseteq\dots\subseteq\M_k$, no $\M_i$ with $i\le k$ contains a witness $Y$ for $b$, so a single constant $\underline{b}$ serves all levels of the tower in $\Gamma$.
A sentence $\forall X\exists Y\varphi$ would instead fail in $\M_k$ at a set $X$ of $\M_k$, which need not belong to the lower models $\M_i$; so there is no counterpart of $\underline{b}$, and the finite subsets of the corresponding $\Gamma$ need not be satisfiable.
In fact, unless $S$ is inconsistent, $S$ is not closed under the rule
\begin{prooftree}
\AxiomC{$\forall Z\forall X\exists Y\varphi(X,Y,Z)$}
\UnaryInfC{$\forall Z\exists\M\bigl(Z\in\M\land\M\models T+\forall X\exists Y\varphi(X,Y,Z)\bigr)$}
\end{prooftree}
for arithmetical $\varphi$.
Indeed, $S$ proves $\wRFN(T)$, which is the premise for $\varphi(X,Y,Z):\equiv X\in Y\land Y\models T$, but an $\omega$-model of $T+\wRFN(T)$ is an $\omega$-model of $S$, so by G\"odel's second incompleteness theorem $S$ does not prove that such a model exists.

Let $\SAC^-$ be $\ACA$ together with the instances $\forall n\exists Y\theta(n,Y)\to\exists W\forall n\theta(n,(W)_n)$ of $\SAC$ for $\Sigma^1_1$-formulas $\theta(n,Y)$ without set parameters.

\begin{corollary}\label{cor:necessary}
Let $\chi$ be a $\Pi^1_2$-sentence such that $\omega$-model reflection is invertible on the cone above $\chi$.
Then $\ACA+\chi\vdash\forall Z\forall n\exists Y\theta(n,Y,Z)\to\forall Z\exists W\forall n\theta(n,(W)_n,Z)$ for every $\Sigma^1_1$-formula $\theta(n,Y,Z)$.
In particular, $\ACA+\chi$ proves $\SAC^-$.
\end{corollary}

\begin{proof}
The $\Pi^1_2$-sentence $\varphi:\equiv\chi\land\forall Z\forall n\exists Y\theta(n,Y,Z)$ lies in the cone above $\chi$, so it is $\ACA$-provably equivalent to $\wRFN(\ACA+\psi)$ for some sentence $\psi$.
Apply Theorem~\ref{thm:B} to $T=\ACA+\psi$; note that $\ACA+\wRFN(T)$ and $\ACA+\varphi$ have the same theorems.
\end{proof}

Let $\mathrm{HYP}$ be the $\omega$-model of the hyperarithmetical sets~\cite[Section~VIII.3]{Sim09}.

\begin{proposition}\label{prop:HYP}
Let $\chi$ be a $\Pi^1_2$-sentence that holds in $\mathrm{HYP}$.
Then there is an $\omega$-model of $\ACA+\chi$ in which an instance of $\SAC^-$ fails.
Hence $\omega$-model reflection is not invertible on the cone above $\chi$.
\end{proposition}

\begin{proof}
The $\omega$-model $\mathrm{HYP}$ satisfies $\SDC$~\cite[Corollary~VIII.4.17]{Sim09}, which proves $\omega$-model reflection for $\Sigma^1_3$-formulas~\cite[Theorem~VIII.5.12]{Sim09}.
So, as in the proof of~\cite[Theorem~5.5]{Has24}, there is a hyperarithmetical set $\M$ that is a coded $\omega$-model of $\ACA+\chi$.
The hyperarithmetical set $\{j:j\notin(\M)_j\}$ is not a set of $\M$.
On the other hand, every $\omega$-model $\M'$ of $\SAC^-$ contains every hyperarithmetical set.
Indeed, every hyperarithmetical set is computable from the jump hierarchy along a computable well-ordering starting with $\varnothing$~\cite[Example~I.11.3]{Sim09}, and $\M'$ contains these hierarchies, by the argument in the proof of Proposition~\ref{lem:cmp}: if $a$ were the least element of a computable well-ordering $\alpha$ such that the hierarchy along the initial segment of $\alpha$ below $a$ is not in $\M'$, then the hierarchies along the initial segments below the elements $b<_\alpha a$ would be in $\M'$; since they are unique, an instance of $\SAC$ without set parameters would collect them in $\M'$, and arithmetical comprehension in $\M'$ would give the hierarchy below $a$.
So some instance of $\SAC^-$ fails in $\M$, and the last claim follows from Corollary~\ref{cor:necessary}.
\end{proof}

For instance, $\ACAp$ and $\Pi^0_\beta\text{-}\mathsf{CA}_0$ for computable well-orderings $\beta\ge1$ hold in $\mathrm{HYP}$, and so does $\forall X\exists Y\theta(X,Y)$ if $\theta$ is arithmetical and every $X$ has a solution hyperarithmetic in $X$.
Together with Theorem~\ref{thm:A}, Proposition~\ref{prop:HYP} also shows that $\ATR$ fails in $\mathrm{HYP}$, which is well known~\cite[Proposition~V.2.6]{Sim09}.
Since $\mathrm{HYP}$ satisfies $\omega$-model reflection for $\Sigma^1_3$-formulas, it satisfies $\wRFN(\ACA+\chi)$ for every $\Pi^1_2$-sentence $\chi$ that holds in it, and none of these sentences implies $\ATR$ over $\ACA$.
For a true sentence $\forall X\exists!Y\theta(X,Y)$ with $\theta$ arithmetical, the unique solution is $\Delta^1_1$ in $X$ and hence hyperarithmetic in $X$, so this recovers~\cite[Corollary~5.13]{Has24}, which says that $\wRFN(\ACA+\forall X\exists!Y\theta(X,Y))$ is not equivalent to $\ATR$ over $\ACA$.

We now give explicit examples of sentences in such cones that are not of the form $\wRFN(\ACA+\psi)$.
Let $F$ be an operation on sets such that the relation $Y=F(X)$ is arithmetical and $X'\le_T F(X)$ for every $X$.
Let $\theta_F(n,H):\equiv(H)_0=\varnothing\land\forall i<n\bigl((H)_{i+1}=F((H)_i)\bigr)$, and let $\mathcal I_F$ be the $\omega$-model whose class of sets is $\{X\subseteq\omega:\exists k\in\omega(X\le_T F^k(\varnothing))\}$.

\begin{proposition}\label{prop:IF}
The $\omega$-model $\mathcal I_F$ satisfies $\ACA+\forall n\exists H\theta_F(n,H)+\neg\exists W\forall n\theta_F(n,(W)_n)$.
Hence, for every $\Pi^1_2$-sentence $\chi$ that holds in $\mathcal I_F$, the sentence $\chi\land\forall n\exists H\theta_F(n,H)$ in the cone above $\chi$ is not $\ACA$-provably equivalent to $\wRFN(\ACA+\psi)$ for any sentence $\psi$, so $\omega$-model reflection is not invertible on the cone above $\chi$.
\end{proposition}

\begin{proof}
Since $F^k(\varnothing)'\le_T F^{k+1}(\varnothing)$ for every $k$, the family $\mathcal I_F$ is closed under $\oplus$, $\le_T$ and the Turing jump, so it is a model of $\ACA$.
Since $\theta_F$ is arithmetical, $\mathcal I_F\models\theta_F(n,H)$ if and only if $\theta_F(n,H)$, for every $H\in\mathcal I_F$.
So $\mathcal I_F$ satisfies $\forall n\exists H\theta_F(n,H)$, because for each $n$ the tuple $\langle F^0(\varnothing),\dots,F^n(\varnothing)\rangle$ is computable from $F^n(\varnothing)$.
It does not satisfy $\exists W\forall n\theta_F(n,(W)_n)$, because such a $W$ would compute $F^n(\varnothing)$ for every $n$, while every set of $\mathcal I_F$ is computable from some $F^k(\varnothing)$, which does not compute $F^{k+1}(\varnothing)$.

Now let $\chi$ be a $\Pi^1_2$-sentence that holds in $\mathcal I_F$, and suppose that $\varphi:\equiv\chi\land\forall n\exists H\theta_F(n,H)$ is $\ACA$-provably equivalent to $\wRFN(\ACA+\psi)$ for some sentence $\psi$.
Since $\ACA+\wRFN(\ACA+\psi)$ proves $\forall n\exists H\theta_F(n,H)$, it proves $\exists W\forall n\theta_F(n,(W)_n)$ by Theorem~\ref{thm:B}, and hence $\ACA+\varphi\vdash\exists W\forall n\theta_F(n,(W)_n)$.
This contradicts the first claim, since $\mathcal I_F\models\varphi$.
\end{proof}

\begin{corollary}\label{cor:examples}
Let $\beta\ge1$ be a computable well-ordering, put $F(X):=X^{(\beta)}$, and let $\chi$ be a $\Pi^1_2$-sentence provable in $\Pi^0_\beta\text{-}\mathsf{CA}_0$.
Then the sentence $\chi\land\forall n\exists H\theta_F(n,H)$ in the cone above $\chi$ is not $\ACA$-provably equivalent to $\wRFN(\ACA+\psi)$ for any sentence $\psi$.
In particular, this applies to $\chi=\ACAp$, and to $\chi=\wRFN(\Pi^1_2\text{-}\omega\mathbf R^\alpha(\ACA))$ for any computable well-ordering $\alpha$.
\end{corollary}

\begin{proof}
The relation $Y=X^{(\beta)}$ is arithmetical, and $X'\le_T X^{(\beta)}$ since $\beta\ge1$.
If $X\le_T F^k(\varnothing)$, then $X^{(\beta)}\le_T F^{k+1}(\varnothing)$, so $\mathcal I_F$ satisfies $\forall X\exists Y(Y=X^{(\beta)})$; since $\beta$ is well-ordered, $\mathcal I_F$ is an $\omega$-model of $\Pi^0_\beta\text{-}\mathsf{CA}_0$, and hence of $\chi$.
So the first claim follows from Proposition~\ref{prop:IF}, and the case of $\ACAp$ is $\beta=\omega$.
For the last case, the sentence $\wRFN(\Pi^1_2\text{-}\omega\mathbf R^\alpha(\ACA))$ is equivalent over $\ACA$ to $\Pi^0_{\omega^{1+\alpha}}\text{-}\mathsf{CA}_0$ by~\cite[Theorem~2.12]{PW23}; apply the first claim with $\beta=\omega^{1+\alpha}$.
\end{proof}

Finally we combine Theorem~\ref{thm:B} with a conservation lemma of Pakhomov and Walsh.
Let $R$ be a $\Pi^1_3$-axiomatized theory, and for $i\in I$ let $\varphi_i$ be a $\Pi^1_2$-formula and $\psi_i$ a $\Pi^1_3$-formula with free variables among $\vec x,\vec X$.
By~\cite[Lemma~2.10]{PW23}, applied with $n=2$ and with all $L_2$-formulas as the countable fragment, if $R$ is closed under the rules
\begin{prooftree}
\AxiomC{$\forall\vec x\forall\vec X\forall\vec y\forall\vec Y(\theta\lor\varphi_i)$}
\UnaryInfC{$\forall\vec x\forall\vec X\forall\vec y\forall\vec Y(\theta\lor\psi_i)$}
\end{prooftree}
for $i\in I$ and $\Pi^1_2$-formulas $\theta(\vec x,\vec X,\vec y,\vec Y)$, then the extension of $R$ by the axioms $\forall\vec x\forall\vec X(\varphi_i\to\psi_i)$ for $i\in I$ is $\Pi^1_2$-conservative over $R$.

\begin{corollary}\label{cor:cons}
Let $T$ and $S$ be as in Theorem~\ref{thm:B}.
Then the extension of $S$ by the sentences $\forall Z\bigl[\forall n\exists Y\varphi(n,Y,Z)\to\exists\M\bigl(Z\in\M\land\M\models T+\forall n\exists Y\varphi(n,Y,Z)\bigr)\bigr]$ for all $\Sigma^1_1$-formulas $\varphi(n,Y,Z)$ is $\Pi^1_2$-conservative over $S$.
In particular, so is $S+\SAC$.
\end{corollary}

\begin{proof}
We apply~\cite[Lemma~2.10]{PW23} as above with $R=S$, where $\varphi_i$ is $\forall n\exists Y\varphi(n,Y,Z)$ and $\psi_i$ is $\exists\M(Z\in\M\land\M\models T+\forall n\exists Y\varphi(n,Y,Z))$.
By pairing we may assume that the side formula $\theta$ has, besides $Z$, exactly one free variable, a set variable $V$.
So we have to show that if $S\vdash\forall Z\forall V(\theta(Z,V)\lor\varphi_i)$, then $S\vdash\forall Z\forall V(\theta(Z,V)\lor\psi_i)$.
Without $\theta$ this is Theorem~\ref{thm:B}, and we modify its proof.
Suppose that $S\vdash\forall Z\forall V(\theta(Z,V)\lor\varphi_i)$ and $S\nvdash\forall Z\forall V(\theta(Z,V)\lor\psi_i)$.
Then there are a model $\U$ of $S$ and sets $Z_0,V_0$ of $\U$ such that $\U\models\neg\theta(Z_0,V_0)$ and, in $\U$, no $\omega$-model of $T$ containing $Z_0$ satisfies $\forall n\exists Y\varphi(n,Y,Z_0)$.
Write $\neg\theta(Z,V)$ as $\exists X\xi(Z,V,X)$ with $\xi$ a $\Pi^1_1$-formula, fix $X_0$ in $\U$ with $\U\models\xi(Z_0,V_0,X_0)$, and let $Z_1$ be a set of $\U$ with $(Z_1)_0=Z_0$, $(Z_1)_1=V_0$ and $(Z_1)_2=X_0$.
We run the proof of Theorem~\ref{thm:B} with $Z_1$ in place of $Z_0$ in the tower and as the interpretation of $\underline{Z}$, with $(\underline{Z})_0$ in place of $\underline{Z}$ in the sentences $\forall j\neg\varphi(\underline{b},(\underline{\M_i})_j,\underline{Z})$, and with the additional sentence $\xi((\underline{Z})_0,(\underline{Z})_1,(\underline{Z})_2)$ in $\Gamma$.
Since $\xi$ is $\Pi^1_1$ and every set of $\N$ is a set of $\U^*$, this sentence holds in $\N$, so $\N\models\neg\theta((\underline{Z})_0,(\underline{Z})_1)$.
As $\N\models S$, we get $\N\models\forall n\exists Y\varphi(n,Y,(\underline{Z})_0)$, and we reach the same contradiction as in the proof of Theorem~\ref{thm:B}.

For the second claim, every instance of $\SAC$ follows over $\ACA$, by pairing, from one of the form $\forall Z(\forall n\exists Y\varphi(n,Y,Z)\to\exists W\forall n\varphi(n,(W)_n,Z))$ with $\varphi$ a $\Sigma^1_1$-formula, and this follows from the corresponding sentence above as in the last paragraph of the proof of Theorem~\ref{thm:B}.
\end{proof}

For $T=\ACA$, Corollary~\ref{cor:cons} says that $\ACAp+\SAC$ is $\Pi^1_2$-conservative over $\ACAp$.
Corollary~\ref{cor:cons} is thus an analogue of the theorem of Barwise and Schlipf~\cite{BS75} that $\SAC$ is $\Pi^1_2$-conservative over $\ACA$ (see also~\cite[Theorem~IX.4.4]{Sim09}).
It does not contain the theorem of Barwise and Schlipf, since $\ACA$ is not of the form $S$.
On the other hand, for a $\Pi^1_2$-sentence $\psi$, the theory $\SAC+\psi$ need not be $\Pi^1_2$-conservative over $\ACA+\psi$.
Indeed, for $F$ as in Proposition~\ref{prop:IF} and $\psi:\equiv\forall n\exists H\theta_F(n,H)$, the theory $\SAC+\psi$ proves the $\Sigma^1_1$-sentence $\exists W\forall n\theta_F(n,(W)_n)$, whereas $\ACA+\psi$ does not by that proposition.
Corollary~\ref{cor:cons} shows that $\SAC+\psi$ is $\Pi^1_2$-conservative over $\ACA+\psi$ for $\psi:\equiv\wRFN(T)$.

\section{\texorpdfstring{The case of $\ATR$}{The case of ATR0}}\label{sec:atr}

Since $\omega$-model reflection is invertible on the cone above $\ATR$ by Theorem~\ref{thm:A}, $\ATR$ itself is of the form $\wRFN(\ACA+\psi)$ for some $\Pi^1_2$-sentence $\psi$.
In this section we give a natural such $\psi$, built from the Veblen functions.

For linear orders $X$ and $Y$ let $\varphi_X(Y)$ be the linear order of the relativized Veblen notation system, written $\varphi(X,Y)$ in~\cite{MM11}; thus, writing $0$ for the empty linear order, $\varphi_0(Y)$ represents $\omega^Y$, and $\varphi_X$ is the $X$-th Veblen function.
The terms of $\varphi_X(Y)$ are built from $0$ and a constant for each element of $Y$ by $+$ and by function symbols $\varphi_x$ for $x\in X$~\cite[Definition~2.7]{MM11}.
Put $\WOP(\varphi_X):\equiv\forall Y(\WO(Y)\to\WO(\varphi_X(Y)))$.
For $x\in X$ we write $X{\restriction}x$ for $\{y\in X:y<x\}$, and $X+1$ for the extension of $X$ by a new largest element.
We write $\varphi_X^n$ for the $n$-fold iterate of $Z\mapsto\varphi_X(Z)$, and for linear orders $X$ and $Y$ we let $\varphi^+_{X+1}(Y):=\sum_{n<\omega}\varphi_X^n(Y+1)$ be the ordered sum of the linear orders $\varphi_X^n(Y+1)$.
If $X$ and $Y$ are well-ordered, the order type of $\varphi^+_{X+1}(Y)$ is the least fixed point of $\varphi_X$ above the order type of $Y$, that is, the least value of $\varphi_{X+1}$ strictly above it.
Pakhomov and Walsh introduce $\varphi^+$ only as a standard notation system for this function~\cite[Section~5.3]{PW23}; we take this choice, and their proofs of~\cite[Theorems~5.12--5.15]{PW23} go through for it.
Following~\cite[Section~5.2]{PW23}, we write $|T|_{\Pi^1_1}$ for the $\Pi^1_1$ proof-theoretic ordinal of a theory $T$ and $|\alpha|$ for the rank of a linear order $\alpha$; equalities and inequalities between these are understood as formalized in $\ACA$ there.
We use the following properties of the notation systems of~\cite{MM11}, which are verified in $\ACA$ by induction on terms.
\begin{enumerate}[label=(V\arabic*),conds]
\item\label{V1} If $K\unlhd L$, then the map on terms induced by an isomorphism of $K$ onto an initial segment of $L$ is an isomorphism of $\varphi_K(0)$ onto an initial segment of $\varphi_L(0)$. Moreover, $x\mapsto\varphi_x(0)$ is an embedding of $X$ into $\varphi_X(0)$.
\item\label{V2} If $X$ is nonempty and has no largest element, then every element of $\varphi_X(0)$ lies in the initial segment $\varphi_{X{\restriction}x}(0)$ for some $x\in X$.
\item\label{V3} For every $n$, $\varphi_X^n(0)\unlhd\varphi_X^{n+1}(0)$, and $\varphi_{X+1}(0)$ embeds into the union of the $\varphi_X^n(0)$.
\item\label{V4} If $e$ is an embedding of $Z$ into $\varphi_{X+1}(0)$, then sending the constant for $z\in Z$ to the term $\varphi_t(e(z))$, where $t$ is the largest element of $X+1$, induces an embedding of $\varphi_X(Z)$ into $\varphi_{X+1}(0)$.
\end{enumerate}
By a theorem of Friedman (see~\cite[Theorem~1.2]{MM11}), $\ACA\vdash\ATR\leftrightarrow\forall X(\WO(X)\to\WO(\varphi_X(0)))$.
Let $\mathsf{I}:\equiv\forall X[\WO(\varphi_X(0))\to\WOP(\varphi_X)]$, which is a $\Pi^1_2$-sentence.

\begin{theorem}\label{thm:atr}
$\ACA\vdash\wRFN(\ACA+\mathsf I)\leftrightarrow\ATR$.
\end{theorem}

\begin{proof}[Proof of the forward direction]
We argue in $\ACA$ and assume $\wRFN(\ACA+\mathsf I)$; by~\cite[Theorem~1.2]{MM11} it suffices to show $\WO(\varphi_X(0))$ for a well-ordering $X$.
Suppose that $D$ is a descending sequence in $\varphi_X(0)$, and let $\M$ be an $\omega$-model of $\ACA+\mathsf I$ containing $X\oplus D$.
Put $K:=X+1$ with largest element $t$, so that $K{\restriction}t=X$.
The set $E$ of $x\in K$ with $\M\models\neg\WO(\varphi_{K{\restriction}x}(0))$ is arithmetical in $\M$, and $t\in E$.
As $K$ is well-ordered, $E$ has a least element $x_0$.
The case where $x_0$ is the least element of $K$ is impossible, since $\varphi_0(0)$ has one element.
If $x_0$ is a limit point of $K$, take a descending sequence in $\varphi_{K{\restriction}x_0}(0)$ belonging to $\M$; by~\ref{V2} its first element lies in the initial segment $\varphi_{K{\restriction}x}(0)$ for some $x<x_0$, hence so does the whole sequence, contradicting the minimality of $x_0$.
If $x_0$ is the successor of $y$, put $L:=K{\restriction}y$, so that $K{\restriction}x_0=L+1$.
Then $\M\models\WO(\varphi_L(0))$, so $\M\models\WOP(\varphi_L)$ by $\mathsf I$.
The set $\{n:\M\models\WO(\varphi_L^n(0))\}$ is arithmetical in $\M$, contains $0$ and is closed under successor; so $\M\models\WO(\varphi_L^n(0))$ for all $n$, and by~\ref{V3}, $\M$ satisfies the well-foundedness of the union of the $\varphi_L^n(0)$.
By~\ref{V3}, $\M\models\WO(\varphi_{L+1}(0))$, that is, $\M\models\WO(\varphi_{K{\restriction}x_0}(0))$, a contradiction.
\end{proof}

For the converse direction we use the results of Pakhomov and Walsh on iterated reflection.
For a $\Pi^1_2$-axiomatized theory $T$, they show~\cite[Theorem~4.3]{PW23} that
\begin{equation}\label{eq:pw-on}
\ACA\vdash\wRFN(T)\leftrightarrow\forall\alpha\bigl(\WO(\alpha)\to\Pi^1_1\text{-}\mathsf{RFN}(\Pi^1_1\text{-}\mathbf R^\alpha(T))\bigr),
\end{equation}
where $\Pi^1_1\text{-}\mathbf R^\alpha(T)$ is $T$ together with $\Pi^1_1\text{-}\mathsf{RFN}(\Pi^1_1\text{-}\mathbf R^\beta(T))$ for all $\beta<\alpha$.
We need five lemmas.

\begin{lemma}\label{lem:mm}
$\ACA\vdash\forall C\bigl[\WO(C)\land\forall Z\exists X\bigl(X=Z^{(\omega^C)}\bigr)\to\WOP(\varphi_C)\bigr]$.
\end{lemma}

\begin{proof}
Marcone and Montalb\'an~\cite[Theorem~3.6]{MM11} show that if $C$ is a $Z$-computable well-ordering, $Y$ is a $Z$-computable linear order and $\varphi_C(Y)$ has a $Z$-computable descending sequence, then $Z^{(\omega^C)}$ computes a descending sequence in $Y$, and they note that the proof formalizes in $\Pi^0_{\omega^\alpha}\text{-}\mathsf{CA}_0$ for each fixed computable well-ordering $\alpha$~\cite[Corollary~3.7]{MM11}.
We check that the proof gives the uniform statement.
We argue in $\ACA$.
Let $\WO(C)$ and $\WO(Y)$, suppose that $D$ is a descending sequence in $\varphi_C(Y)$, and put $Z:=C\oplus Y\oplus D$ and $H:=Z^{(\omega^C)}$.
For $\gamma\in\omega^C$ let $H[<\gamma]$ be the part of $H$ below $\gamma$ as in~\cite[Definition~2.14]{MM11}.
The proof of~\cite[Theorem~3.6]{MM11} defines, by effective transfinite recursion along $C$, a Turing functional which, for $b\in C+1$, takes indices of a linear order $X$ and of a descending sequence in $\varphi_{(C+1){\restriction}b}(X)$, both computable in $H[<\gamma]$, and returns an index of a descending sequence in $X$ computable in $H[<\gamma+\omega^b]$.
For $b\in C+1$ let $\pi(b)$ say that this functional is correct for $b$, for all such indices and all $\gamma$ with $\gamma+\omega^b\le\omega^C$.
Since the sets $H[<\gamma]$ are uniformly computable from $H$, the formula $\pi(b)$ is arithmetical in $H$, and the correctness argument of~\cite{MM11} shows in $\ACA$ that $\pi(b)$ follows from $\pi(b')$ for all $b'<b$.
By arithmetical transfinite induction along the well-ordering $C+1$, $\pi$ holds for the largest element of $C+1$, which concerns $\varphi_C$.
Applying it with $\gamma=0$ to $Y$ and $D$ gives a descending sequence in $Y$, contradicting $\WO(Y)$.
\end{proof}

\begin{lemma}\label{lem:plus}
$\ACA\vdash\forall A\forall C\bigl[A+1\unlhd C\land\WO(\varphi_C(0))\to\WO(\varphi^+_{A+1}(A))\bigr]$.
\end{lemma}

\begin{proof}
We argue in $\ACA$.
By~\ref{V1}, $\varphi_{A+1}(0)$ is isomorphic to an initial segment of $\varphi_C(0)$ and hence well-ordered, and $x\mapsto\varphi_x(0)$ embeds $A+1$ into it.
By~\ref{V4} and induction on $n$, each $\varphi_A^n(A+1)$ embeds into $\varphi_{A+1}(0)$, uniformly in $n$, and so it is well-ordered.
A descending sequence in the ordered sum $\varphi^+_{A+1}(A)$ would eventually stay in one summand, so $\varphi^+_{A+1}(A)$ is well-ordered.
\end{proof}

For a linear order $A$ let $U_A$ be the theory $\ACA+\WO(A)+\forall Z\exists X\bigl(X=Z^{(\omega^A)}\bigr)$.

\begin{lemma}\label{lem:atr-rfn}
$\ACA\vdash\forall A\bigl[\WO(A)\land\WO(\varphi^+_{A+1}(A))\to\Pi^1_1\text{-}\mathsf{RFN}(U_A)\bigr]$.
\end{lemma}

\begin{proof}
We argue in $\ACA$ and assume $\WO(A)$ and $\WO(\varphi^+_{A+1}(A))$.
If $A$ is empty, put $\delta:=0$; then $U_A$ and $\ACA=\Pi^1_2\text{-}\omega\mathbf R^0(\ACA)$ have the same theorems.
Otherwise $A$ has a least element since $\WO(A)$, so $A=1+B$ for a linear order $B$, and we put $\delta:=B+1$.
By~\cite[Theorem~2.12]{PW23}, $\forall Z\exists X(X=Z^{(\omega^A)})$ is $\ACA$-provably equivalent to $\wRFN(\Pi^1_2\text{-}\omega\mathbf R^B(\ACA))$, and by~\cite[Lemma~2.3]{PW23}, adding this sentence to $\ACA$ gives a theory that is equivalent to $\Pi^1_2\text{-}\omega\mathbf R^{B+1}(\ACA)$ up to $\Sigma^1_1$-provability in the sense of~\cite[Section~2.2]{PW23}.
So in both cases $U_A$ is equivalent up to $\Sigma^1_1$-provability to the theory $S:=\Pi^1_2\text{-}\omega\mathbf R^\delta(\ACA)+\WO(A)$, and $1+\delta=A+1$.
Hence $|U_A|_{\Pi^1_1}=|S|_{\Pi^1_1}$ by~\cite[Proposition~5.11]{PW23} and $|S|_{\Pi^1_1}=|\varphi^+_{A+1}(A)|$ by~\cite[Theorem~5.15]{PW23}, and~\cite[Proposition~5.9(1)]{PW23} gives $\Pi^1_1\text{-}\mathsf{RFN}(U_A)$.
\end{proof}

\begin{lemma}\label{lem:atr-incl}
$\ACA\vdash\forall A\mathrm{Pr}_{U_A}\bigl(\gn{\neg\WO(\varphi^+_{A+1}(A))\to\mathsf I\land\WO(\varphi_A(0))}\bigr)$, where $\mathrm{Pr}_{U_A}$ is the usual provability predicate of $U_A$ and $A$ occurs in the code as a set constant.
\end{lemma}

\begin{proof}
We argue in $U_A+\neg\WO(\varphi^+_{A+1}(A))$.
By Lemma~\ref{lem:mm} with $C=A$ we have $\WOP(\varphi_A)$, and in particular $\WO(\varphi_A(0))$.
Now let $C$ satisfy $\WO(\varphi_C(0))$; then $\WO(C)$ by~\ref{V1}.
We compare $A$ and $C$.
The jump hierarchies along $\omega^A$ bound those along $\omega\cdot A$, and these suffice to carry out arithmetical transfinite recursion along $A$ (see~\cite[Section~V.2]{Sim09}).
Defining by such a recursion the partial map $g(a):=\min\bigl(C\setminus g[A{\restriction}a]\bigr)$ as in the proof of~\cite[Theorem~V.6.8]{Sim09}, we obtain $C\unlhd A$ or $A+1\unlhd C$.
If $C\unlhd A$, then $\omega^C$ is an initial segment of $\omega^A$, so $Z^{(\omega^C)}$ exists for every $Z$, and $\WOP(\varphi_C)$ follows from Lemma~\ref{lem:mm}.
If $A+1\unlhd C$, then $\WO(\varphi_C(0))$ gives $\WO(\varphi^+_{A+1}(A))$ by Lemma~\ref{lem:plus}, which is impossible.
Hence $\WO(\varphi_C(0))\to\WOP(\varphi_C)$ for every $C$.
\end{proof}

\begin{lemma}\label{lem:loeb}
Let $U$ be a finitely axiomatized extension of $\ACA$, possibly with set constants, and let $\varphi$ be a sentence such that $\ACA\vdash\varphi\to\Pi^1_1\text{-}\mathsf{RFN}(U)$.
Then $\ACA\vdash\varphi\to\Pi^1_1\text{-}\mathsf{RFN}(U+\neg\varphi)$.
\end{lemma}

\begin{proof}
We first show in $\ACA$ that every $\Pi^1_1$-sentence $\psi$ provable in $U+\neg\varphi$ is provable in $U$.
From the hypothesis, $U$ proves $\varphi\to(\mathrm{Pr}_U(\gn{\psi})\to\psi)$, uniformly in $\psi$, and $U+\neg\varphi\vdash\psi$ gives $U\vdash\neg\varphi\to\psi$.
Hence $U\vdash\mathrm{Pr}_U(\gn{\psi})\to\psi$, and $U\vdash\psi$ by L\"ob's theorem, which is formalizable in $\ACA$.
Now assume $\varphi$; then $\Pi^1_1\text{-}\mathsf{RFN}(U)$ holds, and by what we have shown it yields $\Pi^1_1\text{-}\mathsf{RFN}(U+\neg\varphi)$.
\end{proof}

\begin{proof}[Proof of the converse direction of Theorem~\ref{thm:atr}]
Put $T:=\ACA+\mathsf I$.
By Lemmas~\ref{lem:atr-rfn} and~\ref{lem:loeb} applied to $U=U_A$ and the sentence $\WO(A)\land\WO(\varphi^+_{A+1}(A))$, and by Lemma~\ref{lem:atr-incl}, since $U_A$ proves $\WO(A)$,
\begin{equation}\label{eq:local}
\ACA\vdash\forall A\bigl[\WO(A)\land\WO(\varphi^+_{A+1}(A))\to\Pi^1_1\text{-}\mathsf{RFN}(T+\WO(\varphi_A(0)))\bigr].
\end{equation}
Let $\theta(\alpha)$ be the formula $\WO(\varphi_{\alpha+1}(0))\to\Pi^1_1\text{-}\mathsf{RFN}(\Pi^1_1\text{-}\mathbf R^\alpha(T))$.
We show $\ACA\vdash\forall\alpha\theta(\alpha)$ by L\"ob's theorem over $\ACA$, as in~\cite[Theorems~5.13--5.15]{PW23}.
So we argue in $\ACA$, assume that $\ACA$ proves $\forall\alpha\theta(\alpha)$, and let $\WO(\varphi_{\alpha+1}(0))$.
Let $\psi$ be a $\Pi^1_1$-sentence provable in $\Pi^1_1\text{-}\mathbf R^\alpha(T)$.
The proof uses finitely many axioms $\Pi^1_1\text{-}\mathsf{RFN}(\Pi^1_1\text{-}\mathbf R^{\alpha{\restriction}x_l}(T))$ with $x_1,\dots,x_m\in\alpha$.
If $m=0$, put $A:=0$ and $\xi:\equiv\alpha\text{ is a linear order}$; otherwise let $x^*$ be the $\le_\alpha$-largest of $x_1,\dots,x_m$, and put $A:=(\alpha{\restriction}x^*)+1$ and $\xi:\equiv\alpha\text{ is a linear order}\land\bigwedge_{l=1}^m x_l\le_\alpha x^*$.
Here $\alpha$ occurs in the arithmetical sentence $\xi$ as a set constant, and $\xi$ is true.
By the assumption, each axiom $\Pi^1_1\text{-}\mathsf{RFN}(\Pi^1_1\text{-}\mathbf R^{\alpha{\restriction}x_l}(T))$ is provable in $\ACA+\WO(\varphi_{(\alpha{\restriction}x_l)+1}(0))$, hence, by~\ref{V1}, in $\ACA+\xi+\WO(\varphi_A(0))$.
So the $\Pi^1_1$-sentence $\xi\to\psi$ is provable in $T+\WO(\varphi_A(0))$.
By~\ref{V1}, $\WO(\varphi_{\alpha+1}(0))$ gives $\WO(\alpha+1)$ and hence $\WO(A)$, and since $A+1\unlhd\alpha+1$, it gives $\WO(\varphi^+_{A+1}(A))$ by Lemma~\ref{lem:plus}.
Hence $\xi\to\psi$ is true by~\eqref{eq:local}, and so is $\psi$.
This proves $\forall\alpha\theta(\alpha)$.

Finally, $\ATR$ proves $\WO(\alpha)\to\WO(\varphi_{\alpha+1}(0))$ by~\cite[Theorem~1.2]{MM11}, so together with $\forall\alpha\theta(\alpha)$ it proves the right-hand side of~\eqref{eq:pw-on}, and hence $\wRFN(T)$.
\end{proof}

\section{Concluding remarks}\label{sec:concl}

Besides $\eta_{\ATR}$ and $\mathsf I$, members of $\mathsf{RFN}^{-1}(\ATR)$ include $\SAC$ (see~\cite[Lemma~VIII.4.19]{Sim09} and~\cite[Section~6]{PW23}), $\Delta^1_1\text{-}\mathsf{CA}$, by a theorem of Rathjen~\cite[Theorem~1.7(ii)]{Rat14}, and the theories of~\cite[Proposition~5.8]{Has24}.
The sentences $\eta_{\ATR}$ and $\mathsf I$ also bear on two further points of~\cite{Has24}.
First, Hashimoto~\cite[Section~5]{Has24} states that $\mathsf{RFN}^{-1}(\ATR)$ contains no finitely axiomatizable $\Pi^1_2$-theory; the proof there covers the sentences $\forall X\exists!Y\theta(X,Y)$ of~\cite[Corollary~5.13]{Has24}, while $\eta_{\ATR}$ and $\mathsf I$ are $\Pi^1_2$-sentences in $\mathsf{RFN}^{-1}(\ATR)$.
Second, by~\cite[Theorem~5.5]{Has24}, $\ACA+\mathsf I$ is a natural theory in $\mathsf{RFN}^{-1}(\ATR)$ that is not a theory of hyperarithmetic analysis, which gives an example for~\cite[Question~6.2]{Has24} on the symmetric difference of these two classes.

\begin{question}
If $\omega$-model reflection is invertible on the cone above a $\Pi^1_2$-sentence $\chi$, does $\chi$ imply $\ATR$ over $\ACA$?
\end{question}

A positive answer would give $\ACA+\chi\vdash\SAC$, since $\ATR$ proves $\SAC$, whereas Corollary~\ref{cor:necessary} gives only the implications $\forall Z\forall n\exists Y\theta(n,Y,Z)\to\forall Z\exists W\forall n\theta(n,(W)_n,Z)$, and in particular $\SAC^-$.
Even this weaker conclusion is open: an instance of $\SAC$ with a set parameter is a $\Pi^1_3$-sentence, so it is not in the cone above $\chi$, and Theorem~\ref{thm:B} does not apply to it.
By Proposition~\ref{prop:HYP}, such a $\chi$ fails in $\mathrm{HYP}$, as $\ATR$ does.

The proof of Theorem~\ref{thm:A} uses $\ATR$ only through Proposition~\ref{lem:cmp}, and the formula proved there implies $\ATR$ over $\ACA$.
Indeed, if $K$ and $L$ are well-orderings and $\M$ is as there for $K$ and $Z=L$, then $\M\models\WO(L)$, so $\M\models L\unlhd K\lor K\lhd L$ and $K$ and $L$ are comparable, and the comparability of well-orderings implies $\ATR$~\cite[Theorem~V.6.8]{Sim09}.
So such a $\chi$ not implying $\ATR$ would require a different argument.

The question of the introduction is trivial over $\SDC$: this theory is equivalent over $\ACA$ to $\Pi^1_3\text{-}\omega\mathsf{RFN}(\ACA)$ (see~\cite[Theorem~VIII.5.12]{Sim09}), so it proves $\varphi\leftrightarrow\wRFN(\ACA+\varphi)$ for every $\Pi^1_2$-sentence $\varphi$.
However, no $\Pi^1_2$-sentence $\chi$ that is consistent with $\ACA$ implies $\SDC$ over $\ACA$: otherwise $\ACA+\chi$ would prove $\wRFN(\ACA+\chi)$ and hence its own consistency, since $\omega$-models have full satisfaction classes in $\ACAp$.
Thus, if $\ACA+\chi$ is consistent and $\omega$-model reflection is invertible on the cone above $\chi$, then $\ACA+\chi$ proves $\SAC^-$ by Corollary~\ref{cor:necessary} but does not prove $\SDC$, as is the case for $\chi=\ATR$.

\section*{Acknowledgements}

The author thanks Yudai Suzuki and Keita Yokoyama for helpful discussions, and the Mathematisches Forschungsinstitut Oberwolfach for its hospitality during the workshop \emph{Mathematical Logic: Proof Theory, Constructive Mathematics} (Workshop 2647, 15--20 November 2026).

\section*{Use of AI}

The results of this paper were obtained through discussions between the author and large language models, mainly Claude Opus~5.5 (Anthropic) and GPT-6.2 Astra (OpenAI); Claude Fable~5.2 (also Anthropic) suggested Proposition~\ref{prop:HYP}, and the manuscript was prepared with Claude Opus~5.5 and Claude Fable~5.2.
The author takes full responsibility for the content of this paper.

\end{document}